\documentclass[12pt]{article}
\usepackage{amsmath,amssymb,amsthm}

\newtheorem{thm}{Theorem}[section] 

\newtheorem{conj}[thm]{Conjecture}

\begin{document}

\date{\small Mathematics Subject Classification: 11B65, 33B99}

\title{Beyond Zudilin's Conjectured $q$-analog of Schmidt's problem}

\author{Thotsaporn ``Aek'' Thanatipanonda\\
\small{\tt thotsaporn@gmail.com}%
}

\maketitle

\begin{abstract}
Using the methodology of (rigorous) {\it experimental mathematics},
we give a simple and motivated
solution to Zudilin's question concerning a $q$-analog
of a problem posed by Asmus Schmidt about a certain binomial coefficients sum.
Our method is based on two
simple identities that can be automatically proved using
the Zeilberger and $q$-Zeilberger algorithms. We further illustrate our method by
proving two further binomial coefficients sums.

\end{abstract}

\section{Introduction}
In 1992, Asmus Schimdt \cite{Schmidt1}, conjectured that
for any integer $r \geq 1$, the sequence of numbers
$\{{c^{(r)}_{k}}\}_{k \geq 0} $
defined implicitly by
\[
\sum_k\binom{n}{k}^r\binom{n+k}{k}^r =
\sum_k\binom{n}{k}\binom{n+k}{k}c^{(r)}_k,\quad n=0,1,2,\dots
\]
are always integers.
It took more than 10 years to completely solve
this conjecture(\cite{Zudilin}).
Shorter proofs were found recently (\cite{Guo,TT}).
The key step in \cite{TT} was the fact that
$a^{(r)}_{k,j}$,  defined implicitly by

\begin{equation} \label{schimdt}
\binom{n}{k}^r\binom{n+k}{k}^r =
\sum_j a^{(r)}_{k,j}\binom{n}{j}\binom{n+j}{j},
\quad n=0,1,2,\dots
\end{equation}

are all integers.\\

The integral property of the $a^{(r)}_{k,j}$,
in fact, directly solves the conjecture.
It is very possible that had Schimdt
conjectured this equation in the first place,
his conjecture would have been proved much sooner.  \\

Inspired by (\ref{schimdt}),
we investigate binomial coefficients summands
of the form \[\prod_i\binom{n+d_ik}{b_ik+c_i},\]
where $d_i,b_i,c_i$ are fixed integers,
such that their powers  can be written
as an integer linear combination of themselves.

\section{Results}
In this section, we search for $f(n,k)$
for which $a^{(r)}_{k,j}$

\begin{equation} \label{main}
f(n,k)^r = \sum_j a^{(r)}_{k,j}f(n,j),\quad n=0,1,2,\dots
\end{equation}

are all integers.\\

It was not an accident that the $a^{(r)}_{k,j}$ happen to be
integers. Experiments show that there exists
an integer-valued function $S_f(k,j,i)$,
free of $r$, such that $\bar{a}^{(r)}_{k,j}$ define by
$ \bar{a}^{(1)}_{k,k}=1$,
$ \bar{a}^{(1)}_{k,j}=0$ for $j\neq k$, and

\begin{equation}
\bar{a}^{(r+1)}_{k,j}
= \sum_{i} S_f(k,j,i) \bar{a}^{(r)}_{k,i}.
\end{equation}

agrees with $a^{(r)}_{k,j}$ in (\ref{main}). \\

We prove this fact by showing that $\bar{a}^{(r)}_{k,j}$
also satisfy (\ref{main}). The proof relies on an induction.
On the one hand,

\begin{align*}
 \sum_{j}\bar{a}^{(r+1)}_{k,j}f(n,j)
          & = \sum_{j}\sum_{i}
                S_f(k,j,i)\bar{a}^{(r)}_{k,i}f(n,j)
          \kern10pt\text{\small (by definition of $\bar{a}^{(r+1)}_{k,j}$)} \\
          & = \sum_{i}\bar{a}^{(r)}_{k,i}\sum_{j} S_f(k,j,i)f(n,j).
\end{align*}

On the other hand,
\begin{align*}
 f(n,k)^{r+1} &  = f(n,k)^{r}f(n,k)\\
          & = \sum_{i}\bar{a}^{(r)}_{k,i} f(n,i)f(n,k)
          \kern10pt\text{\small (by induction hypothesis).}
\end{align*}

Hence the proof  boils down to the condition that $S_f(k,j,i)$
must satisfy:

\begin{equation} \label{con}
f(n,i)f(n,k) = \sum_{j} S_f(k,j,i)f(n,j), \;\ n = 0,1,2,3,... .
\end{equation}

Let's state this observation as a theorem.\\

\begin{thm} Given a pair $f(n,k)$ and $S(k,j,i)$
such that
\[
f(n,i)f(n,k) = \sum_{j} S(k,j,i)f(n,j), \;\ \mbox{for all} \;\ n,i,k \geq 0.
\]

Define $a^{(r)}_{k,j}$ recursively by
$ a^{(1)}_{k,k}=1$, $a^{(1)}_{k,j}=0$ for $j\neq k$ and
\[
\displaystyle a^{(r+1)}_{k,j}
= \sum_{i} S(k,j,i)a^{(r)}_{k,i}.
\]

Then for $k \geq 0$ and $r \geq 1$,
\[
 f(n,k)^r
= \sum_j  a^{(r)}_{k,j}
\displaystyle f(n,j), \;\ n = 0,1,2,... .
\]
\end{thm}

Once we know where to look, of course in this case guided by
Schidmt's problem, it becomes a routine job that computers are
so good at. We first pick some binomial term $f(n,k)$. Then
we crank out some numerical values of $S(k,j,i)$ according to (\ref{con}).
Once we have enough data, we ask our computer
to guess the relation, or even the
formula for $S(k,j,i)$. Finally, needless to say,
the proof of the identities can be routinely done by
Zeilberger's algorithm. Here is the list that we were able to find.\\

\textbf{Result 2.1.1:}\\
These are the functions we used in Schmidt's conjecture.\\
\[
f(n,k) = \binom{n}{k}\binom{n+k}{k}, \;\
S(k,j,i) = \binom{i+k}{i}\binom{j}{j-i,j-k,i+k-j}.
\]

\textbf{Result 2.1.2:}\\
For any fixed integer $c$,
\[
f(n,k) = \binom{n}{k+c}\binom{n+k}{k+c}, \;\
S(k,j,i) = \frac{(i+k+c)!}{(i+c)!(k+c)!}
\frac{(j+c)!}{j!}\binom{j+c}{j-i,j-k,i+k+c-j}.
\]

\textbf{Result 2.2:}\\
\[
f(n,k) = \binom{n}{k}, \;\
S(k,j,i) = \binom{j}{j-i,j-k,i+k-j}.
\]

\textbf{Result 2.3:}\\
\[
f(n,k) = \binom{n+k}{k}, \;\
S(k,j,i) = (-1)^{i+j+k}\binom{j}{j-i,j-k,i+k-j}.
\]

\section{$q$-analog}
    We now  present $q$-analogs of
    the results from section 2.\\

The $q$-binomial $\displaystyle {n \brack k}$ are defined by

\[ {n \brack k} =  \left\{ \begin{array}{ll}
                                \frac{(q)_n}{(q)_k(q)_{n-k}}
                                & \mbox{if} \;\ 0 \leq k \leq n.\\
                                0 & \mbox{otherwise}
                        \end{array} \right.
\]

where $(q)_0=1$ and
$(q)_n=(1-q)(1-q^2)...(1-q^n)$ for $n=1,2,...$.\\

The proof of Theorem 3.1 below is similar to the
proof of Theorem 2.1. We leave the details to the reader.\\

\begin{thm} Given a triple $f(n,k)$, $A(k,j,i,n)$ and $S(k,j,i)$
satisfying
\begin{equation} \label{qcon}
f(n,i)f(n,k) = \sum_{j} q^A S(k,j,i)f(n,j), \;\ \mbox{for all} \;\ n,i,k \geq 0.
\end{equation}
Let $B(k,j,i)$ and $C(k,j,r,n)$ be any functions such that \\

$B(k,j,i)+C(k,j,r+1,n)=A(k,j,i,n)+C(k,i,r,n)$ and  \\

$C(k,k,1,n)=0.$ \\

Define $P^{(r)}_{k,j}(q)$ recursively by
$ P^{(1)}_{k,k}(q)=1$, $P^{(1)}_{k,j}(q)=0$ for $j\neq k$ and
\[
\displaystyle P^{(r+1)}_{k,j}(q)
= \sum_{i} q^B S(k,j,i)P^{(r)}_{k,i}(q).
\]

Then for $k \geq 0$ and $r \geq 1$,
\[
 f(n,k)^r
= \sum_j q^C
\displaystyle f(n,j) P^{(r)}_{k,j}(q), \;\ n = 0,1,2,... .
\]
\end{thm}

The proofs of (\ref{qcon}) of the results below
can again be done automatically using the $q$-Zeilberger
algorithm. Once we find function $A$ from (\ref{qcon}),
it is only a matter of simple
calculations to solve for functions $B$ and $C$.  \\

\textbf{Result 3.1.1:}
\[f(n,k) = {n \brack k}{n+k \brack k}, \;\
S(k,j,i) = {i+k \brack i}{j \brack j-i,j-k,i+k-j},
\]

$ A(k,j,i,n)= (n-j)(k+i-j),$

$B(k,j,i) \;\ \;\ = -(k+i-j)j,$

$C(k,j,r,n) = (rk-j)n.$ \\

\textbf{Result 3.1.2:} \\
\[f(n,k) = {n \brack k+c}{n+k \brack k+c}, \;\
S(k,j,i) = \frac{(q)_{i+k+c}}{(q)_{i+c}(q)_{k+c}}
\frac{(q)_{j+c}}{(q)_{j}}{j \brack j-i,j-k,i+k-j},
\]

$ A(k,j,i,n) = (n-j)(k+i-j)+c(n-k-i-c),$

$ B(k,j,i) \;\ \;\ = -(k+i-j)j-c(k+i+c),$

$ C(k,j,r,n) = (rk-j)n+rcn-cn.$ \\

\textbf{Result 3.2:} \\
\[
f(n,k) = {n \brack k}, \;\
S(k,j,i) = {j \brack j-i,j-k,i+k-j},
\]

$ A(k,j,i,n) = (j-i)(j-k),$

$ B(k,j,i) \;\ \;\ = (j-i)(j-k),$

$ C(k,j,r,n) = 0.$ \\

\textbf{Result 3.3:} \\
\[
f(n,k) = {n+k \brack k}, \;\
S(k,j,i) = (-1)^{i+j+k}{j \brack j-i,j-k,i+k-j},
\]

$ A(k,j,i,n) = \frac{(k+i-j)(2n+k+i-j+1)}{2},$

$ B(k,j,i) \;\ \;\ = \frac{(k+i-j)(k+i-j+1)}{2},$

$ C(k,j,r,n) = (rk-j)n.$ \\

Note that 3.1.1 agrees with the results in \cite{Guo}.
Also, the positivity of $B$ and $C$ in result 3.2 and 3.3 imply that
${P}^{(r)}_{k,j}$ are indeed  polynomials.

\section{Conjectures}
After some experimentation with all of the
binomial terms of the form $\prod_i\binom{n+d_ik}{b_ik+c_i}$,
we found out that the only such terms that satisfy the
condition (\ref{con}) seem to be of the form $\binom{n}{k+c}\binom{n+k}{k+c}$
or $\binom{n+dk}{k}$ for any fixed integer $c$ and $d$.
We already saw that it was true for the former, but
for the latter, it is only proved to be true for
$d = 0,1$, but we conjecture that it holds for all non-negative integers $d$ as follows.\\

\begin{conj}
For any integers $d$,
$k \geq 0$ and $r \geq 1$,
there exist integers ${a}^{(r)}_{d,k,j}$ such that
\[
\binom{n+dk}{k}^r = \sum_j {a}^{(r)}_{d,k,j}\binom{n+dj}{j}
\;\ \;\ \text{for all $n=0,1,2,...$.}
\]
Moreover ${a}^{(r)}_{d,k,j}$ can be defined as following:
$ {a}^{(1)}_{d,k,k}=1$, $ {a}^{(1)}_{d,k,j}=0$ for $j\neq k$ and
\[
{a}^{(r+1)}_{d,k,j} = \sum_{i} S_d(k,j,i) {a}^{(r)}_{d,k,i},
\]
where, $S_d(k,j,i)$ are integers, independent of $r$, for all $d,k,j,i$.
\end{conj}

\begin{conj} For a fixed integer $d$,
$S_d(k,j,i)$ defined above are holonomic
but not of the first order,
ie. no closed form solution, except $d = 0,1$.
\end{conj}

\section{Conclusions}
We presented a motivated and streamlined new proof of
the main result of \cite{Guo}, as well as two new, much deeper, identities,
and made two conjectures.
But the main interest of this paper is in illustrating a {\it methodology},
using computers (via {\it experimental mathematics}), to generate data, then
formulate conjectures, and finally having the very same computer
{\it rigorously} prove its own conjectures. We believe that this
methodology has great potential almost everywhere in mathematics.


\end{document}